\documentclass{article}
\usepackage{amsmath,amssymb,amsthm}
\usepackage{bm}

\newtheorem{theorem}{Theorem}[section]
\newtheorem{proposition}[theorem]{Proposition}
\newtheorem{lemma}[theorem]{Lemma}
\newtheorem{corollary}[theorem]{Corollary} 
\theoremstyle{definition}
\newtheorem{definition}[theorem]{Definition} 
\theoremstyle{remark}

\newcommand{\N}{\mathbb{N}}
\newcommand{\Z}{\mathbb{Z}}

\newcommand{\C}{\mathbb{C}}

\newcommand{\bmu}{\bm{\mu}}

\newcommand{\bx}{\bm{x}}
\newcommand{\bone}{\bm{1}}

\newcommand{\Xvec}{X_1,\ldots,X_r}
\newcommand{\Yvec}{Y_1,\ldots,Y_r}
\newcommand{\nvec}{n_1,\ldots,n_r}
\newcommand{\kvec}{k_1,\ldots,k_r}
\newcommand{\tvec}{t_1,\ldots,t_{p-1}}
\newcommand{\xvec}{x_1,\ldots,x_r}
\newcommand{\nfracvec}{n_1!\, \cdots n_r!}
\newcommand{\kfracvec}{k_1!\, \cdots k_r!}
\newcommand{\zerovec}{0,\ldots,0}
\newcommand{\op}{\partial_{X_i} + \partial_{Y_i} - 1}

\begin{document}
\title{A generalization of the duality for multiple harmonic sums}
\author{Gaku Kawashima}
\date{}
\maketitle

\begin{abstract}
 The duality is a fundamental property of the finite multiple harmonic 
 sums (MHS).
 In this paper, we prove a duality result for certain generalizations of
 MHS which appear naturally as the differences of MHS. 
 We also prove a formula for the differences of these generalized MHS.
\end{abstract}

\begin{center}
 Keywords: multiple harmonic sum\\
% 2000 Mathematics Subject Classification: Primary 11B83 
\end{center}

\section{Introduction}
In recent years, multiple harmonic sums (MHS for short) have been
studied by theoretical physicists and the duality for MHS
was discovered~\cite[Section 8B]{V}.
The same formula has appeared in~\cite{H,Z} in the study of
arithmetical properties of MHS.
It is also used in~\cite{K} in order to investigate algebraic relations among
multiple zeta values.
\par
Here, we explain the duality for MHS.
We call an ordered set of positive integers a multi-index.
If we introduce the following nested sums
\begin{displaymath}
 s_{\bmu}(n) = \sum_{n=n_1 \ge \cdots \ge n_p \ge 0}
 \frac{1}{(n_1+1)^{\mu_1} \cdots (n_p+1)^{\mu_p}}, \quad
 0 \le n \in \Z
\end{displaymath}
for a multi-index $\bmu = (\mu_1, \ldots, \mu_p)$, 
the duality for MHS is expressed as follows:
\begin{equation}
 \sum_{k=0}^{n} (-1)^k \binom{n}{k} s_{\bmu}(k) = s_{\bmu^{*}}(n), \quad
 0 \le n \in \Z, \label{eq0-10}
\end{equation}
where $\bmu^{*}$ is some multi-index determined by $\bmu$.
For example, we have
\begin{displaymath}
  (1,2,3)^{*} = (2,2,1,1),\quad (2,2,2)^{*} = (1,2,2,1)\quad \text{and}\quad (4,1,1)^{*} = (1,1,1,3)
\end{displaymath}
by the diagrams
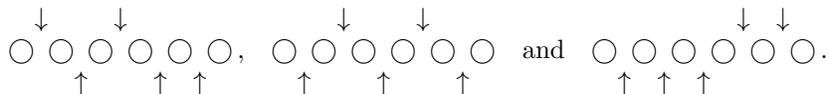
\begin{displaymath}
 {\setlength{\arraycolsep}{0pt}
 \begin{array}{cccccccccccccccccccccccccccccccccccc}
  &\downarrow& & & &\downarrow& & & & & & & & & &\downarrow& & & &\downarrow& & & & & & & & & & & &\downarrow& &\downarrow& &
   \\
  \bigcirc& &\bigcirc& &\bigcirc& &\bigcirc& &\bigcirc& &\bigcirc&\, ,\quad&\bigcirc& &\bigcirc& &\bigcirc& &\bigcirc& &\bigcirc& &\bigcirc&\quad\text{and}\quad&\bigcirc& &\bigcirc& &\bigcirc& &\bigcirc& &\bigcirc& &\bigcirc&\,. \\
  & & &\uparrow& & & &\uparrow& &\uparrow& & & &\uparrow& & & &\uparrow& & & &\uparrow& & & &\uparrow& &\uparrow& &\uparrow& & & & & & \\
 \end{array}}
\end{displaymath}
The lower arrows are in the complementary slots to the upper arrows
(see~\cite[Section 2]{K} for the precise definition of the
correspondence $\bmu \mapsto \bmu^{*}$).
A q-analogue of the identity (\ref{eq0-10}) is also known~\cite{B}.
\par
In~\cite{KT}, the finite nested sums
\begin{displaymath}
 c_{x_1,\ldots,x_p}(n) = \sum_{n=n_1 \ge \cdots \ge n_p \ge 0}
 \frac{x_1^{n_1-n_2} \cdots x_{p-1}^{n_{p-1}-n_p} x_p^{n_p}}
 {(n_1+1) \cdots (n_{p-1}+1)}, \quad 0 \le n \in \Z 
\end{displaymath}
with complex parameters $x_1,\ldots,x_p$ were studied
in order to derive relations among multiple $L$-values
and a generalization of (\ref{eq0-10})
\begin{equation}
 \sum_{k=0}^{n} (-1)^k \binom{n}{k} c_{x_1,\ldots,x_p}(k) 
  = c_{1-x_1,\ldots,1-x_p}(n), \quad 0 \le n \in \Z \label{eq0-20}
\end{equation}
was proved.
To obtain the claim (\ref{eq0-10}) from (\ref{eq0-20}), 
we have only to note that
\begin{align*}
 &c_{0,\ldots,0,1,\ldots,0,\ldots,0,1,0}(n) \quad \text{and} \quad
 c_{0,\ldots,0,1,\ldots,0,\ldots,0,1,0,\ldots,0,1}(n) \\[-4mm]
 &\hspace{1.1mm}\underbrace{\hspace{9mm}}_{\mu_1}
 \hspace{3.9mm}\underbrace{\hspace{9mm}}_{\mu_p}
 \hspace{20.8mm}\underbrace{\hspace{9mm}}_{\mu_1}
 \hspace{3.9mm}\underbrace{\hspace{9mm}}_{\mu_{p-1}}
 \hspace{0.8mm}\underbrace{\hspace{7mm}}_{\mu_{p}}
\end{align*} 
are both equal to $s_{\bmu}(n)$.
The duality (\ref{eq0-20}) is a consequence of the difference formula
\begin{displaymath}
 (\Delta^k c_{x_1,\ldots,x_p})(n)
 = c_{x_1,\ldots,x_p; 1-x_1,\ldots,1-x_p}(n,k), \quad 0 \le n,\,k \in \Z.
\end{displaymath}
The right-hand side is given by the following nested sums
\begin{multline*}
 c_{x_1,\ldots,x_p; y_1,\ldots,y_p}(n,k) = 
 \sum_{\substack{n=n_1 \ge \cdots \ge n_p \ge 0\\
 k=k_1 \ge \cdots \ge k_p \ge 0}} P \\
 \times 
 \frac{(x_1^{n_1-n_2} \cdots x_{p-1}^{n_{p-1}-n_p} x_p^{n_p})
 (y_1^{k_1-k_2} \cdots y_{p-1}^{k_{p-1}-k_p} y_p^{k_p})}
 {(n_1+k_1+1) \cdots (n_{p-1}+k_{p-1}+1)}, \quad 
 0 \le n,\,k \in \Z
\end{multline*}
with complex parameters $x_i$ and $y_i$ for $1 \le i \le p$, where
\begin{displaymath}
 P = \binom{n+k}{n}^{-1}
 \binom{n_1-n_2+k_1-k_2}{n_1-n_2} \cdots 
 \binom{n_{p-1}-n_p+k_{p-1}-k_p}{n_{p-1}-n_p}
 \binom{n_p+k_p}{n_p}.
\end{displaymath}
The symbol $\Delta$ denotes the difference operator 
and its definition is given in Section \ref{sec1}.
\par
In the present paper, we ask whether the duality exists for these more 
complicated nested sums which appear naturally as the differences of
$c_{x_1,\ldots,x_p}(n)$.
We consider more general nested sums
\begin{multline*}
  c_{\bx_1;\cdots;\bx_r}^{\tvec}(\nvec)
  = \sum_{\substack{n_1=n_{11} \ge \cdots \ge n_{1p} \ge 0\\ \cdots\\
  n_r=n_{r1} \ge \cdots \ge n_{rp} \ge 0}} Q \\
  \times
  \frac{(x_{11}^{n_{11}-n_{12}} \cdots x_{1p-1}^{n_{1p-1}-n_{1p}}
  x_{1p}^{n_{1p}}) \cdots
  (x_{r1}^{n_{r1}-n_{r2}} \cdots x_{rp-1}^{n_{rp-1}-n_{rp}} 
  x_{rp}^{n_{rp}})}
  {(n_{11} + \cdots + n_{r1} + t_1) \cdots 
  (n_{1p-1} + \cdots + n_{rp-1} + t_{p-1})}, \\
  0 \le n_1,\, \ldots,\, n_r \in \Z
 \end{multline*}
with parameters 
$t_1$, \ldots, $t_{p-1} \in \C \setminus \{0,-1,-2,\ldots \}$ and
$\bx_i = (x_{i1},\ldots,x_{ip}) \in \C^p$ for $1 \le i \le r$.
We need to find an appropriate $Q$ so that the duality holds.
In addition, it is necessary for $Q$ to become $P$ on setting
$r=2$ and $t_1 = \cdots = t_{p-1} = 1$.
The following
\begin{multline*}
  Q =
  \binom{n_{1}}{n_{11}-n_{12}, \ldots, n_{1p-1}-n_{1p}, n_{1p}} \cdots
  \binom{n_{r}}{n_{r1}-n_{r2}, \ldots, n_{rp-1}-n_{rp}, n_{rp}} \\
  \times
  \binom{n_{11} + \cdots + n_{r1} + t_1 - 1 }
  {n_{11}-n_{12} + \cdots + n_{r1}-n_{r2}}^{-1} \cdots 
  \binom{n_{1p-1} + \cdots + n_{rp-1} + t_{p-1} - 1}
  {n_{1p-1}-n_{1p} + \cdots + n_{rp-1}-n_{rp}}^{-1},
\end{multline*}
which consists of multinomial coefficients and generalized binomial
coefficients,
satisfies such requirements.
In fact, for the above $Q$ we have
\begin{multline*}
 \sum_{k_1=0}^{n_1} \cdots \sum_{k_r=0}^{n_r}
 (-1)^{k_1+\cdots+k_r} \binom{n_1}{k_1} \cdots \binom{n_r}{k_r}
 c_{\bx_1;\cdots;\bx_r}^{\tvec}(\kvec) \\
 = c_{\bone-\bx_1;\cdots;\bone-\bx_r}^{\tvec}(\nvec), \quad
 0 \le n_1,\, \ldots,\, n_r \in \Z,
\end{multline*}
where $\bone-\bx := (1-x_1,\ldots,1-x_p)$ for any
$\bx = (x_1,\ldots,x_p) \in \C^p$. 
This is the main result Corollary \ref{cor2-70} of this paper.
Moreover, we shall determine the differences of
$c_{\bx_1;\cdots;\bx_r}^{\tvec}(\nvec)$ in Corollary \ref{cor3-30}
explicitly.
\section{The differences of multiple sequences} \label{sec1}
In the present section, we study some fundamental properties 
of the differences of multiple sequences. 
For this purpose, we use formal power series belonging to
$\C[[X]] := \C[[X_1,\ldots,X_r]]$ or 
$\C[[X,Y]] := \C[[X_1,\ldots,X_r,Y_1,\ldots,Y_r]]$, where $r$ is a
positive integer and is fixed throughout this section.
We put $\partial_{X_i}=\partial/\partial X_i$ and
$\partial_{Y_i}=\partial/\partial Y_i$ for $1 \le i \le r$.
In the following, we denote the set of non-negative integers by $\N$.
We sometimes write $n$ (resp. $k$) for a sequence $n_1,\ldots,n_r$
(resp. $k_1,\ldots,k_r$) for abbreviation.
First, we note that for a formal power series
\begin{displaymath}
 F = \sum_{n,k=0}^{\infty} a(\nvec,\kvec)
 \frac{X_1^{n_1} \cdots X_r^{n_r} Y_1^{k_1} \cdots Y_r^{k_r}}{\nfracvec \,
 \kfracvec} 
 \in \C[[X,Y]],
\end{displaymath}
we have the equality
\begin{multline}
 (\op)F =\sum_{n,k=0}^{\infty} 
 \Bigl\{a(n_1,\ldots,n_i+1,\ldots,n_r,k)
 + a(n, k_1,\ldots,k_i+1,\ldots,k_r) \Bigr.\\
 \Bigl. - a(n,k) \Bigr\} 
 \frac{X_1^{n_1} \cdots X_r^{n_r} Y_1^{k_1} \cdots Y_r^{k_r}}{\nfracvec \,
 \kfracvec} \label{eq1-10}
\end{multline}
for any $1 \le i \le r$.
\begin{lemma}
 \label{lem1-20}
 If a formal power series $F\in \C[[X,Y]]$ satisfies two conditions
 \begin{displaymath}
  (\op) F = 0 \quad (\text{for any } 1 \le i \le r)
 \end{displaymath}
 and
 \begin{displaymath}
  F(\Xvec,\zerovec) = 0,
 \end{displaymath}
 then we have $F = 0$.
\end{lemma}
\begin{proof}
 If we put
 \begin{displaymath}
  F = \sum_{n,k=0}^{\infty} a(\nvec,\kvec)
  \frac{X_1^{n_1} \cdots X_r^{n_r} Y_1^{k_1} \cdots Y_r^{k_r}}{\nfracvec \,
  \kfracvec}, 
 \end{displaymath}
 by (\ref{eq1-10}) we have
 \begin{displaymath}
  a(n_1,\ldots,n_i+1,\ldots,n_r,k) + a(n,k_1,\ldots,k_i+1,\ldots,k_r)
  - a(n,k) =0
 \end{displaymath}
 for any $\nvec,\kvec \in \N$ and any $1 \le i \le r$. Moreover we have
 \begin{displaymath}
  a(\nvec,\zerovec) = 0
 \end{displaymath}
 for any $\nvec \in \N$. By induction on $k_1 + \cdots + k_r$, we see that
 \begin{displaymath}
  a(\nvec,\kvec) = 0
 \end{displaymath}
 holds for any $\nvec,\kvec \in \N$.
\end{proof}
\begin{definition}
 \label{df1-30}
 We denote by $\C^{\N^r}$ the set of all mappings from $\N^r$ to $\C$.
 For each $1 \le i \le r$, we define the difference operator
 $\Delta_i \colon \C^{\N^r} \to \C^{\N^r}$ by putting
 \begin{displaymath}
  (\Delta_i a)(\nvec) = a(\nvec) - a(n_1,\ldots,n_i+1,\ldots,n_r)
 \end{displaymath}
 for any $a \in \C^{\N^r}$ and any $\nvec \in \N$.
\end{definition}
We note that $\Delta_i \Delta_j = \Delta_j \Delta_i$ for any 
$1 \le i,\,j \le r$. We denote the composition
\begin{displaymath}
 \underbrace{\Delta_1 \circ \cdots \circ \Delta_1}_{n_1} \circ \underbrace{\Delta_2 \circ \cdots \circ
 \Delta_2}_{n_2} \circ \cdots \circ \underbrace{\Delta_r \circ \cdots \circ \Delta_r}_{n_r}
\end{displaymath}
by $\Delta_1^{n_1} \Delta_2^{n_2} \cdots \Delta_r^{n_r}$.
\begin{definition}
 \label{df1-40}
 We define the inversion operator
 $\nabla \colon \C^{\N^r} \to \C^{\N^r}$ by putting
 \begin{displaymath}
  (\nabla a)(\nvec) = (\Delta_1^{n_1} \cdots \Delta_r^{n_r} a)(\zerovec)
 \end{displaymath}
 for any $a \in \C^{\N^r}$ and any $\nvec \in \N$.
\end{definition}
For a sequence $a \colon \N^r \to \C$, we define
\begin{displaymath}
 f_a = \sum_{n=0}^{\infty} a(\nvec) \frac{X_1^{n_1} \cdots
 X_r^{n_r}}{\nfracvec} \in \C[[X]]
\end{displaymath}
and
\begin{displaymath}
 F_a = \sum_{n,k=0}^{\infty} (\Delta_1^{k_1} \cdots
 \Delta_r^{k_r} a)(\nvec) \frac{X_1^{n_1} \cdots X_r^{n_r} Y_1^{k_1}
 \cdots Y_r^{k_r}}{\nfracvec \, \kfracvec} 
 \in \C[[X,Y]].
\end{displaymath}
We note that
\begin{equation}
 F_a(\Xvec,\zerovec) = f_a(\Xvec) \label{eq1-50}
\end{equation}
and
\begin{equation}
 F_a(\zerovec,\Yvec) = f_{\nabla a}(\Yvec). \label{eq1-60}
\end{equation}
By (\ref{eq1-10}), we have
\begin{equation}
 (\op)F_a = 0 \label{eq1-70}
\end{equation}
for any $1 \le i \le r$.
\begin{proposition}
 \label{prop1-80}
 Let $a \colon \N^r \to \C$ be a sequence. Then we
 have
 \begin{displaymath}
  F_{\nabla a}(\Xvec,\Yvec) = F_a(\Yvec,\Xvec).
 \end{displaymath}
\end{proposition}
\begin{proof}
 By (\ref{eq1-70}), we have
 \begin{displaymath}
  (\op)F_{\nabla a}(\Xvec,\Yvec) = 0
 \end{displaymath}
and
\begin{displaymath}
 (\op)F_a(\Yvec,\Xvec) = 0
\end{displaymath}
for any $1 \le i \le r$. In addition, by (\ref{eq1-50}) and
 (\ref{eq1-60}) we have
\begin{displaymath}
 F_{\nabla a}(\Xvec,\zerovec) = f_{\nabla a}(\Xvec)
\end{displaymath}
and
\begin{displaymath}
 F_a(\zerovec,\Xvec) = f_{\nabla a}(\Xvec).
\end{displaymath}
Therefore Lemma \ref{lem1-20} implies the proposition.
\end{proof}
\begin{corollary}
 \label{cor1-90}
 Let $a \colon \N^r \to \C$ be a sequence. Then for any
 $\nvec$, $\kvec \in \N$, we have
 \begin{displaymath}
  (\Delta_1^{k_1} \cdots \Delta_r^{k_r}(\nabla a))(\nvec) =
  (\Delta_1^{n_1} \cdots \Delta_r^{n_r} a)(\kvec).
 \end{displaymath}
\end{corollary}
\begin{proof}
 By comparing the coefficients of both sides of the equation in
 Proposition \ref{prop1-80}, we obtain the corollary.
\end{proof}
\begin{corollary}
 \label{cor1-100}
 The operator $\nabla^2$ is the identity on $\C^{\N^r}$.
\end{corollary}
\begin{proof}
 It follows from Corollary \ref{cor1-90} on setting $n_1=\cdots = n_r = 0$.
\end{proof}
In the rest of this section, we give explicit expressions for the
differences and the inversions of multiple sequences.
\begin{proposition}
 \label{prop1-110}
 For any sequence $a \colon \N^r \to \C$, we have
 \begin{equation*}
  F_a = f_a(X_1-Y_1,\ldots,X_r-Y_r) e^{Y_1+\cdots +Y_r}. 
 \end{equation*}
\end{proposition}
\begin{proof}
 It is easily seen that for any $1 \le i \le r$ the right-hand side
 becomes zero if we apply $\op$.
 Noting (\ref{eq1-50}) and (\ref{eq1-70}),
 we get the proposition by Lemma \ref{lem1-20}.
\end{proof}
\begin{corollary}
 \label{cor1-120}
 For any sequence $a \colon \N^r \to \C$ and any
 $\nvec$, $\kvec \in \N$, we have
 \begin{multline*}
  (\Delta_1^{k_1} \cdots \Delta_r^{k_r} a)(\nvec) \\
  = \sum_{i_1=0}^{k_1} \cdots \sum_{i_r=0}^{k_r} (-1)^{i_1 + \cdots +
  i_r} \binom{k_1}{i_1} \cdots \binom{k_r}{i_r} a(n_1+i_1,\ldots,n_r+i_r).
 \end{multline*}
\begin{proof}
 If we apply the operator $\partial_{Y_1}^{k_1} \cdots \partial_{Y_r}^{k_r}$
 to the right-hand side of the equation in Proposition \ref{prop1-110}, 
 the result is
 \begin{displaymath}
  \sum_{i_1=0}^{k_1} \cdots \sum_{i_r=0}^{k_r}
  (-1)^{i_1+\cdots+i_r} \binom{k_1}{i_1} \cdots \binom{k_r}{i_r}
  f_a^{(i_1,\ldots,i_r)}(X_1-Y_1,\ldots,X_r-Y_r) e^{Y_1 + \cdots + Y_r}
 \end{displaymath}
 by the Leibniz rule, where we have put
 \begin{displaymath}
  f_a^{(i_1,\ldots,i_r)} = \partial_{X_1}^{i_1} \cdots
  \partial_{X_r}^{i_r} f_a.
 \end{displaymath}
 Therefore by applying the operator
 $\partial_{X_1}^{n_1} \cdots \partial_{X_r}^{n_r} \partial_{Y_1}^{k_1} \cdots \partial_{Y_r}^{k_r}$
 to both sides of the equation in Proposition \ref{prop1-110}
 and comparing the constant term, we obtain the desired equality.
\end{proof}
\end{corollary}
\begin{corollary}
 \label{cor1-140}
 Let $a \colon \N^r \to \C$ be a sequence. For any $\nvec \in \N$, we have
 \begin{displaymath}
  (\nabla a)(\nvec) 
  = \sum_{i_1=0}^{n_1} \cdots \sum_{i_r=0}^{n_r}
  (-1)^{i_1+\cdots+i_r} \binom{n_1}{i_1} \cdots \binom{n_r}{i_r}
  a(i_1,\ldots,i_r).
 \end{displaymath}
\end{corollary}
\begin{proof}
 It follows from Corollary \ref{cor1-120} on setting 
 $n_1 = \cdots = n_r = 0$.
\end{proof}
\section{A proof of the duality} \label{sec2}
Throughout this section, we fix positive integers $r$ and $p$.
First of all, we define and calculate the inversions of operators on
$\C[[X]] := \C[[\Xvec]]$.
Then, we shall prove a generalization of the duality for MHS 
by using this notion.
We begin with the definition of the inversions of formal power series.
For a formal power series
\begin{displaymath}
 f = \sum_{n=0}^{\infty}a(\nvec)
 \frac{X_1^{n_1} \cdots X_r^{n_r}}{n_1! \cdots n_r!} \in \C[[X]],
\end{displaymath}
we define its inversion by
\begin{displaymath}
 \nabla f = \sum_{n=0}^{\infty}(\nabla a)(\nvec)
 \frac{X_1^{n_1} \cdots X_r^{n_r}}{n_1! \cdots n_r!}.
\end{displaymath}
By Corollary \ref{cor1-100}, the inversion operator $\nabla$ on 
$\C[[X]]$ is an involution.
\begin{proposition}
 \label{prop2-10}
 For any formal power series $f \in \C[[X]]$, we have
 \begin{displaymath}
  \nabla f = f(-X_1,\ldots,-X_r) e^{X_1+\cdots+X_r}.
 \end{displaymath}
\end{proposition}
\begin{proof}
 It is immediate from (\ref{eq1-60}) and Proposition \ref{prop1-110}.
\end{proof}
For a mapping $\xi \colon \C[[X]] \to \C[[X]]$, we define its inversion
by
\begin{displaymath}
 (\nabla \xi)f = \nabla(\xi(\nabla f)), \quad f \in \C[[X]].
\end{displaymath}
Clearly, we have
\begin{displaymath}
 \nabla(\xi f) = (\nabla \xi)(\nabla f)
\end{displaymath}
for any $\xi \colon \C[[X]] \to \C[[X]]$ and any $f \in \C[[X]]$
since the inversion operator $\nabla$ on $\C[[X]]$ is an involution.
It is easily seen that the equality
\begin{equation}
 \nabla(\xi\eta) = (\nabla\xi)(\nabla\eta) \label{eq2-10}
\end{equation}
holds for any mappings $\xi$, $\eta \colon \C[[X]] \to \C[[X]]$.
We note that the inversion operator $\nabla$ is a linear operator.
\begin{proposition}
 \label{prop2-20}
 For any $1 \le i \le r$, we have
 \begin{displaymath}
  \nabla X_i = -X_i \quad \text{and} \quad 
  \nabla \partial_{X_i} = 1 - \partial_{X_i},
 \end{displaymath}
 where $X_i$ is the operator on $\C[[X]]$ defined by
 $f \mapsto X_i f$ for any $f \in \C[[X]]$.
\end{proposition}
\begin{proof}
 By Proposition \ref{prop2-10}, for any $f \in \C[[X]]$ we have
 \begin{displaymath}
  (\nabla X_i)(\nabla f) = \nabla(X_i f)
  = -X_i f(-X_1,\ldots,-X_r) e^{X_1+\cdots+X_r}
  = -X_i (\nabla f),
 \end{displaymath}
 which implies $\nabla X_i = -X_i$.
 If we calculate $(\nabla \partial_{X_i})f$ by definition,
 the other assertion also follows easily.
\end{proof}
In the following, we consider a generalization of the duality for MHS.
We recall that we have fixed positive integers $r$ and $p$.
Let $\bx_i = (x_{i1},\ldots,x_{ip}) \in \C^p$ for $1 \le i \le r$ and let
$\tvec \in \C \setminus \{0,-1,-2,\ldots\}$.
Then for any $\nvec \in \N$, we define 
$c_{\bx_1;\cdots;\bx_r}^{\tvec}(\nvec)$ by
\begin{displaymath}
  \sum_{\substack{n_1=n_{11} \ge \cdots \ge n_{1p} \ge 0\\ \cdots\\
  n_r=n_{r1} \ge \cdots \ge n_{rp} \ge 0}} \:
  \frac{\displaystyle{\prod_{i=1}^r
  \binom{n_{i}}{\nu_{i1},\ldots,\nu_{ip}}
  x_{i1}^{\nu_{i1}} \cdots x_{ip}^{\nu_{ip}}}}
  {\displaystyle{\prod_{j=1}^{p-1} 
  \binom{n_{1j}+\cdots+n_{rj}+t_j-1}{\nu_{1j}+\cdots+\nu_{rj}}
  (n_{1j} + \cdots + n_{rj} + t_j)}},
\end{displaymath}
where we have put
\begin{displaymath}
 \nu_{i1} = n_{i1} - n_{i2}, \quad \ldots, \quad
 \nu_{ip-1} = n_{ip-1} - n_{ip}, \quad \nu_{ip} = n_{ip}
\end{displaymath}
for abbreviation.
In the case $p=1$, we have
\begin{equation}
 c_{x_1;\cdots;x_r}(\nvec) = x_1^{n_1} \cdots x_r^{n_r}. \label{eq2-20}
\end{equation}
The following Proposition \ref{prop2-30} is the key 
to proving the duality for the
nested sums $c_{\bx_1;\cdots;\bx_r}^{\tvec}(\nvec)$.
For any $\bx = (x_1,\ldots,x_p) \in \C^p$ with $p \ge 2$, we define
\begin{displaymath}
 {}^{-}\bx = (x_2,\ldots,x_p) \in \C^{p-1}.
\end{displaymath}
\begin{proposition}
 \label{prop2-30}
 Let $\bx_i = (x_{i1},\ldots,x_{ip}) \in \C^p$ for $1 \le i \le r$ and
 let $\tvec \in \C \setminus \{0,-1,-2,\ldots\}$. If $p \ge 2$, for 
 any $\nvec \in \N$ we have
 \begin{multline}
  (n_1 + \cdots + n_r + t_1) 
  c_{\bx_1;\cdots;\bx_r}^{t_1,\ldots,t_{p-1}}(\nvec) \\
  - \sum_{k=1}^r x_{k1} n_k
  c_{\bx_1;\cdots;\bx_r}^{\tvec}(n_1,\ldots,n_k-1,\ldots,n_r)
  = c_{{}^{-}\bx_1;\cdots;{}^{-}\bx_r}^{t_2,\ldots,t_{p-1}}(\nvec).
  \label{eq2-25}
 \end{multline}
\end{proposition}
\begin{proof}
 Let $1 \le k \le r$. 
 In the proof we write $\sum_{\sharp}$ and $\sum_{\flat}$ instead of
 \begin{displaymath}
  \sum_{\substack{n_1 \ge n_{12} \ge \cdots \ge n_{1p} \ge
  0\\ \cdots\\ n_k > n_{k2} \ge \cdots \ge n_{kp} \ge 0\\ \cdots \\
  n_r \ge n_{r2} \ge \cdots \ge n_{rp} \ge 0}}
  \quad \text{and} \quad 
  \sum_{\substack{n_1 \ge n_{12} \ge \cdots \ge n_{1p} \ge 0\\ \cdots\\
  n_r \ge n_{r2} \ge \cdots \ge n_{rp} \ge 0\\ \exists i,\, n_i > n_{i2}}},
 \end{displaymath}
 respectively. 
 The $k$th term in the second sum of the left-hand side of (\ref{eq2-25})
 is equal to
 \begin{equation}
  \sum_{\sharp} R_k
  \times
  \frac{\displaystyle{x_{k1}^{\nu_{k1}} \cdots x_{kp}^{\nu_{kp}}
  \prod_{\substack{1 \le i \le r \\ i \ne k}}
  \binom{n_{i}}{\nu_{i1},\ldots,\nu_{ip}}
  x_{i1}^{\nu_{i1}} \cdots x_{ip}^{\nu_{ip}}}}
  {\displaystyle{
  \prod_{j=2}^{p-1} 
  \binom{n_{1j}+\cdots+n_{rj}+t_j-1}{\nu_{1j}+\cdots+\nu_{rj}}
  (n_{1j} + \cdots + n_{rj} + t_j)}}, \label{eq2-30}
 \end{equation}
 where
 \begin{displaymath}
  R_k = n_k \binom{n_k-1}{\nu_{k1}-1,\nu_{k2},\ldots,\nu_{kp}}
  \binom{n_1+\cdots+n_r+t_1-2}{\nu_{11}+\cdots+\nu_{r1}-1}^{-1}
  \frac{1}{n_1+\cdots+n_r+t_1-1}.
 \end{displaymath}
 By a simple calculation we obtain
 \begin{displaymath}
  R_k = \frac{\nu_{k1}}{\nu_{11}+\cdots+\nu_{r1}}
  \binom{n_k}{\nu_{k1},\ldots,\nu_{kp}}
  \binom{n_1+\cdots+n_r+t_1-1}{\nu_{11}+\cdots+\nu_{r1}}^{-1},
 \end{displaymath}
 which is obviously equal to $0$ if $n_k=n_{k2}$ (i.e. $\nu_{k1}=0$).
 Therefore we can replace $\sharp$ by $\flat$ in the expression (\ref{eq2-30}).
 Consequently the second sum of the left-hand side of (\ref{eq2-25}) 
 is equal to
 \begin{displaymath}
  \sum_{\flat}
  \frac{\displaystyle{
  \prod_{i=1}^r
  \binom{n_{i}}{\nu_{i1},\ldots,\nu_{ip}}
  x_{i1}^{\nu_{i1}} \cdots x_{ip}^{\nu_{ip}}}}
  {\displaystyle{
  \binom{n_1+\cdots+n_r+t_1-1}{\nu_{11}+\cdots+\nu_{r1}}
  \prod_{j=2}^{p-1} 
  \binom{n_{1j}+\cdots+n_{rj}+t_j-1}{\nu_{1j}+\cdots+\nu_{rj}}
  (n_{1j} + \cdots + n_{rj} + t_j)}},
 \end{displaymath}
 from which the desired equality is easily derived.
\end{proof}
In order to rewrite Proposition \ref{prop2-30} in terms of
generating functions, we introduce the following operator on $\C[[X]]$
\begin{displaymath}
 \xi_{\xvec} = X_1 \partial_{X_1} + \cdots + X_r \partial_{X_r} -
 x_1 X_1 - \cdots - x_r X_r
\end{displaymath}
for complex numbers $\xvec$.
For this operator we have
\begin{equation}
 \nabla \xi_{\xvec} = \xi_{1-x_1,\ldots,1-x_r} \label{eq2-35}
\end{equation}
by Proposition \ref{prop2-20}.
We note that for a formal power series
\begin{displaymath}
 f = \sum_{n=0}^{\infty} a(\nvec) \frac{X_1^{n_1} \cdots
 X_r^{n_r}}{\nfracvec} \in \C[[X]]
\end{displaymath}
and a complex number $t$, the equality
\begin{multline}
 (\xi_{\xvec} + t) f
 = \sum_{n=0}^{\infty} \Bigl\{(n_1+\cdots+n_r+t) a(\nvec)
 \Bigr. \\
 \Bigl. - \sum_{k=1}^r x_k n_k a(n_1,\ldots,n_k-1,\ldots,n_r)\Bigr\}
 \frac{X_1^{n_1} \cdots X_r^{n_r}}{\nfracvec} \label{eq2-50}
\end{multline}
holds.
\begin{lemma}
 \label{lem2-40}
 For any $\xvec \in \C$ and any $t \in \C \setminus \{0,-1,-2,\ldots\}$, the
 operator $\xi_{\xvec} + t$ on $\C[[X]]$ is an injection.
\end{lemma}
\begin{proof}
 Let
 \begin{displaymath}
  f = \sum_{n=0}^{\infty} a(\nvec) \frac{X_1^{n_1} \cdots
  X_r^{n_r}}{\nfracvec} \in \C[[X]]
 \end{displaymath}
 be a formal power series in the kernel of the operator 
 $\xi_{\xvec} + t$.
 Then, by (\ref{eq2-50}), the equality
 \begin{displaymath}
  (n_1+\cdots+n_r+t) a(\nvec) = \sum_{k=1}^r x_k n_k
  a(n_1,\ldots,n_k-1,\ldots,n_r) 
 \end{displaymath}
 holds for any $\nvec \in \N$. Since we have
 \begin{displaymath}
  n_1+\cdots+n_r+t \neq 0
 \end{displaymath}
 for any $\nvec \in \N$, we obtain $f=0$
 by induction on $n_1 + \cdots + n_r$.
\end{proof}
For any $\bx_1,\ldots,\bx_r \in \C^p$ and any
$\tvec \in \C \setminus \{0,-1,-2,\ldots\}$, we put
\begin{displaymath}
 f_{\bx_1;\cdots;\bx_r}^{\tvec}
 = \sum_{n=0}^{\infty} c_{\bx_1;\cdots;\bx_r}^{\tvec}(\nvec)
 \frac{X_1^{n_1} \cdots X_r^{n_r}}{\nfracvec} 
 \in \C[[X]].
\end{displaymath}
In the case $p=1$, we have
\begin{equation}
 f_{x_1;\cdots;x_r} = e^{x_1 X_1 + \cdots + x_r X_r} \label{eq2-60}
\end{equation}
by (\ref{eq2-20}).
\begin{proposition}
 \label{prop2-50}
 Let $\bx_i=(x_{i1},\ldots,x_{ip}) \in \C^p$ for $1 \le i \le r$ and let
 $\tvec \in \C \setminus \{0,-1,-2,\ldots\}$.
 If $p \ge 2$ then we have
 \begin{displaymath}
  (\xi_{x_{11},\ldots,x_{r1}} + t_1)
  f_{\bx_1;\cdots;\bx_r}^{\tvec} =
  f_{{}^{-}\bx_1;\cdots;{}^{-}\bx_r}^{t_2,\ldots,t_{p-1}}.
 \end{displaymath}
\end{proposition}
\begin{proof}
 It follows from Proposition \ref{prop2-30} and (\ref{eq2-50}).
\end{proof}
Now we state the main result of this paper.
We put
\begin{displaymath}
 \bone-\bx = (1-x_1,\ldots,1-x_p)
\end{displaymath}
for any $\bx = (x_1,\ldots,x_p) \in \C^p$.
\begin{theorem}
 \label{th2-60}
 For any $\bx_1,\ldots,\bx_r \in \C^p$ and any 
 $\tvec \in \C \setminus \{0,-1,-2,\ldots\}$, we have
 \begin{displaymath}
  \nabla f_{\bx_1;\cdots;\bx_r}^{\tvec}
  = f_{\bone-\bx_1;\cdots;\bone-\bx_r}^{\tvec}.
 \end{displaymath}
\end{theorem}
\begin{proof}
 By Proposition \ref{prop2-50} and (\ref{eq2-60}), we have
 \begin{displaymath}
  (\xi_{x_{1p-1},\ldots,x_{rp-1}} + t_{p-1}) \cdots
  (\xi_{x_{11},\ldots,x_{r1}} + t_{1})
  f_{\bx_1;\cdots;\bx_r}^{\tvec} = e^{x_{1p}X_1 + \cdots + x_{rp}X_r}.
 \end{displaymath}
 Applying the inversion operator $\nabla$ to both sides, we obtain
 \begin{multline*}
  (\xi_{1-x_{1p-1},\ldots,1-x_{rp-1}} + t_{p-1}) \cdots
  (\xi_{1-x_{11},\ldots,1-x_{r1}} + t_{1})
  (\nabla f_{\bx_1;\cdots;\bx_r}^{\tvec}) \\
  = e^{(1-x_{1p})X_1 + \cdots + (1-x_{rp})X_r}
 \end{multline*}
 by (\ref{eq2-10}), (\ref{eq2-35}) and Proposition \ref{prop2-10}.
 On the other hand we have
 \begin{multline*}
  (\xi_{1-x_{1p-1},\ldots,1-x_{rp-1}} + t_{p-1}) \cdots
  (\xi_{1-x_{11},\ldots,1-x_{r1}} + t_{1})
  f_{\bone-\bx_1;\cdots;\bone-\bx_r}^{\tvec} \\
  = e^{(1-x_{1p})X_1 + \cdots + (1-x_{rp})X_r},
 \end{multline*} 
 which together with Lemma \ref{lem2-40} implies the theorem.
\end{proof}
\begin{corollary}
 \label{cor2-70}
 For any $\bx_1,\ldots,\bx_r \in \C^p$ and any
 $\tvec \in \C \setminus \{0,-1,-2,\ldots\}$, we have
 \begin{displaymath}
  \nabla c_{\bx_1;\cdots;\bx_r}^{\tvec} 
  = c_{\bone-\bx_1;\cdots;\bone-\bx_r}^{\tvec}.
 \end{displaymath}
\end{corollary}
\begin{proof}
 It is clear from Theorem \ref{th2-60}.
\end{proof}
\section{The difference formula} \label{sec3}
In this section, we shall prove a formula 
which gives explicitly the differences of
the nested sums $c_{\bx_1;\cdots;\bx_r}^{\tvec}(\nvec)$.
As usual, we put
\begin{displaymath}
 [\xi,\eta] = \xi \eta - \eta \xi
\end{displaymath}
for mappings $\xi$, $\eta \colon \C[[X]] \to \C[[X]]$.
\begin{lemma}
 \label{lem3-10}
 Let $i_1, \ldots, i_q$ be distinct integers satisfying 
 $1 \le i_1,\ldots,i_q \le r$
 and let $\xvec,t$ be complex numbers.
 If we put
 \begin{displaymath}
  x_{i_1} + \cdots + x_{i_q} = c,
 \end{displaymath}
 then we have
 \begin{displaymath}
  [\partial_{X_{i_1}} + \cdots + \partial_{X_{i_q}} - c,\, \xi_{\xvec} +
  t] = \partial_{X_{i_1}} + \cdots + \partial_{X_{i_q}} - c.
 \end{displaymath}
\end{lemma}
\begin{proof}
 The left-hand side is equal to
 \begin{displaymath}
  [\partial_{X_{i_1}},\, X_{i_1}\partial_{X_{i_1}}] + \cdots +
  [\partial_{X_{i_q}},\, X_{i_q}\partial_{X_{i_q}}] -
  x_{i_1}[\partial_{X_{i_1}},\, X_{i_1}] - \cdots -
  x_{i_q}[\partial_{X_{i_q}},\, X_{i_q}].
 \end{displaymath}
 Since we have
 $[\partial_{X_j},\, X_j] = 1$ and
 $[\partial_{X_j},\, X_j \partial_{X_j}] = \partial_{X_j}$
 for all $1 \le j \le r$, we obtain the assertion.
\end{proof}
\begin{theorem}
 \label{th3-20}
 Let $\bx_1,\ldots,\bx_r \in \C^p$ and let $c \in \C$.
 We suppose that
 \begin{displaymath}
  \bx_{i_1} + \cdots + \bx_{i_q} = \underbrace{(c,\ldots,c)}_p
 \end{displaymath}
 for distinct integers $1 \le i_1,\ldots,i_q \le r$. 
 Then for any $\tvec \in \C \setminus \{0,-1,-2,\ldots\}$, we have
 \begin{displaymath}
  (\partial_{X_{i_1}} + \cdots + \partial_{X_{i_q}} - c)
  f_{\bx_1;\cdots;\bx_r}^{\tvec} = 0.
 \end{displaymath}
\end{theorem}
\begin{proof}
 The proof is by induction on $p$. 
 The case $p=1$ follows directly from (\ref{eq2-60}).
 Let $p \ge 2$. 
 We put $\bx_i = (x_{i1},\ldots,x_{ip})$ for any $1 \le i \le r$.
 Since we have
 \begin{displaymath}
  (\xi_{x_{11},\ldots,x_{r1}} + t_1 + 1)(\partial_{X_{i_1}} + \cdots +
  \partial_{X_{i_q}} - c)
  =
  (\partial_{X_{i_1}} + \cdots + \partial_{X_{i_q}} -
  c)(\xi_{x_{11},\ldots,x_{r1}} + t_1)
 \end{displaymath}
 by Lemma \ref{lem3-10}, it holds that
 \begin{displaymath}
  (\xi_{x_{11},\ldots,x_{r1}} + t_1 + 1)(\partial_{X_{i_1}} + \cdots +
  \partial_{X_{i_q}} - c)f_{\bx_1;\cdots;\bx_r}^{\tvec} = 0
 \end{displaymath}
 by Proposition \ref{prop2-50} and the hypothesis of induction.
 This, together with Lemma \ref{lem2-40}, leads to
 \begin{displaymath}
  (\partial_{X_{i_1}} + \cdots + \partial_{X_{i_q}} - c)
  f_{\bx_1;\cdots;\bx_r}^{\tvec} = 0.
 \end{displaymath}
 Therefore we complete the proof.
\end{proof}
\begin{corollary}
 \label{cor3-30}
 Let $\bx_1,\ldots,\bx_r \in \C^p$
 and let $\tvec \in \C \setminus \{0,-1,-2,\ldots\}$. Then, for
 any $\nvec,\kvec \in \N$, we have
 \begin{displaymath}
  (\Delta_1^{k_1} \cdots \Delta_r^{k_r}
  c_{\bx_1;\cdots;\bx_r}^{\tvec})(\nvec) = c_{\bx_1;\cdots;\bx_r;\bone-\bx_1;\cdots;\bone-\bx_r}^{\tvec}(\nvec,\kvec).
 \end{displaymath}
\end{corollary}
\begin{proof}
 By Theorem \ref{th3-20}, we have
 \begin{displaymath}
  (\partial_{X_i} + \partial_{Y_i} - 1)
  f_{\bx_1;\cdots;\bx_r;\bone-\bx_1;\cdots;\bone-\bx_r}^{\tvec}(\Xvec,\Yvec)
  = 0
 \end{displaymath}
 for any $1 \le i \le r$. We also have
 \begin{displaymath}
  f_{\bx_1;\cdots;\bx_r;\bone-\bx_1;\cdots;\bone-\bx_r}^{\tvec}(\Xvec,\zerovec)
  = f_{\bx_1;\cdots;\bx_r}^{\tvec}(\Xvec),
 \end{displaymath}
 which results from
 \begin{displaymath}
  c_{\bx_1;\cdots;\bx_r;\bone-\bx_1;\cdots;\bone-\bx_r}^{\tvec}(\nvec,\zerovec)
  = c_{\bx_1;\cdots;\bx_r}^{\tvec}(\nvec).
 \end{displaymath}
 Therefore by Lemma \ref{lem1-20}, (\ref{eq1-50}) and (\ref{eq1-70}) 
 we obtain the assertion.
\end{proof}
By this corollary, we see that 
the differences of the nested sums considered in this paper
are again the same kind of nested sums.
Corollary \ref{cor2-70} follows also from Corollary \ref{cor3-30}
on setting $n_1=\cdots=n_r=0$.

\begin{flushleft}
 Graduate School of Mathematics \\
 Nagoya University \\
 Chikusa-ku, Nagoya 464-8602, Japan \\
 E-mail: m02009c@math.nagoya-u.ac.jp
\end{flushleft}

\begin{thebibliography}{9}
  \bibitem{B}
    D. M. Bradley,
    Duality for finite multiple harmonic q-series,
    Disc. Math. {\bf{300}} (2005), 44-56.
  \bibitem{H}
   M. Hoffman, Quasi-symmetric functions and mod $p$ multiple harmonic sums,
   preprint arXiv:math.NT/0401319.
  \bibitem{K}
   G. Kawashima, A class of relations among multiple zeta values,
   J. Number Theory. {\bf{129}} (2009), 755-788.
  \bibitem{KT}
   G. Kawashima and T. Tanaka, Newton series and extended derivation
   relations for multiple $L$-values,
   preprint arXiv:0801.3062.
  \bibitem{V}
    J. A. M. Vermaseren,
    Harmonic sums, Mellin transforms and integrals,
    Int. J. Mod. Phys. A {\bf{14}} (1999), 2037-2076.
  \bibitem{Z}
    J. Zhao,
    Wolstenholme type theorem for multiple harmonic sums,
    Int. J. Number Theory {\bf{4}} (2008), 73-106.	
\end{thebibliography}
\end{document}